\documentclass[12pt]{article}
\usepackage[cp1251]{inputenc}
\usepackage[russian]{babel}
\usepackage{amsfonts}

\newtheorem{lemma}{Lemma}
\newtheorem{theorem}{Theorem}

\begin{document}

\def\C{{\mathbb C}}
\def\R{{\mathbb R}}
\def\N{{\mathbb N}}
\def\i{{\sqrt {-1}}}
\def\O{{\cal O}}
\def\A{{\cal A}}
\def\F{{\cal F}}
\def\M{{\widehat{M}}}
\def\DW{{\widehat{\cal D}}}
\def\D{{\cal D}}
\def\Z{{\mathbb Z}}
\def\L{{\cal L}}
\def\A{{\cal A}}
\def\B{{\cal B}}
\def\P{{\cal P}}
\def\Pic{\rm{Pic}}
\def\Mat{\rm{Mat}}
\def\ord{\rm{ord}}
\def\Im{\rm{Im}}

\title{On new examples of Hamiltonian--minimal and minimal
Lagrangian submanifolds in ${\mathbb C}^n$ and ${\mathbb CP}^n$ }
\author{A. E. Mironov} \date{}

\maketitle

\section{Introduction and main results}

In this article we propose a new method for the construction of
Hamiltonian--minimal (H-minimal) and minimal Lagrangian
immersions  of some manifolds in $\C^n$ and in $\C P^n$. By this method one can construct, in particular,
immersions of such manifolds as the generalized Klein's bottle ${\cal K}^n$,
the multi\-di\-men\-sional torus, ${\cal K}^{n-1}\times S^1$, $S^{n-1}\times S^1$,
and others. In some cases these immersions are embeddings. For
example, it is possible to embed the following manifolds:  ${\cal K}^{2n+1},$
$S^{2n+1}\times S^1$, ${\cal K}^{2n+1}\times S^1$,
$S^{2n+1}\times S^1\times S^1$.

An immersion $\psi : L \rightarrow P$ of an
$n$-dimensional manifold $L$ in a symplectic $2n$-dimensional
manifold $P$ with a symplectic form $\omega$ is
called Lagrangian, if $\psi^*(\omega)=0.$

In the sequel let $P$ be a K\"ahler manifold with a symplectic
form
$\omega=-{\rm Im}ds^2,$ where $ds^2$ is a Hermitian
metric on $P$.

A Lagrangian immersion $\psi$ is called
$H$-minimal if the variations of volume $\psi(L)$ along all
Hamiltonian fields  with compact supports vanish:
$$
\frac{d}{d t}{\rm vol}(\psi_t(L))\vert_{t=0}=0,\eqno{(1)}
$$
where
$\psi_0(L)=\psi(L), \psi_t(L)$ is a deformation
$\psi(L)$ along a Hamiltonian field $W$, i.e. such a field that the
one-form
$\omega_W(\psi^*(\xi))=\psi^*\omega(W,\xi)$ is exact on
$L$, and ${\rm vol}(\psi_t(L))$ is the volume of the deformed part of
$\psi_t(L)$.

As an example of Hamiltonian deformations we can consider
deformations of a unit circle in a unit sphere
$S^2$ when $S^2$ is laid out on two parts of equal area in
each moment of time. We recall that an immersion is called
minimal if the variations  of volume $\psi(L)$ along all fields
vanish.

Let $M$ be a $k$-dimensional manifold given in
$\R^n$ by a system of equations
$$
e_{1j}u_1^2+\dots+e_{nj}u_n^2=d_j,\ j=1,\dots,n-k, \eqno{(2)}
$$
where $d_j\in\R,e_{ij}\in\Z.$

Since ${\rm dim}M=k$, we can assume that equations (2)
are linearly independent and the integer vectors
$e_j=(e_{j1},\dots,e_{j(n-k)})\in\Z^{n-k},j=1,\dots,n$
form a lattice
$\Lambda$
in
$\R^{n-k}$ of maximum rank. Let's denote by
$\Lambda^*$ a dual lattice to $\Lambda$, i.e.,
$$
\Lambda^*=\{\lambda^*\in\R^{n-k}|
(\lambda^*,\lambda)\in\Z,\lambda\in\Lambda\}.
$$
By $\Gamma$ we denote the following factor group
$$
 \Gamma=\Lambda^*/2\Lambda^*.
$$
 The following isomorphism takes place
$$
 \Gamma = \Z_2^{n-k}.
$$
By $T^{n-k}$ we denote the
$(n-k)$-dimensional torus
$$
  T^{n-k}=\{(e^{\pi i(e_1,y)},\dots,e^{\pi i(e_n,y)})\}\subset \C^n,
$$
 where
$y=(y_1,\dots,y_{n-k}) \in \R^{n-k}$ and
$(e_j,y)=e_{j1}y_1+\dots+e_{j(n-k)}y_{n-k}$.

We define an action of group $\Gamma$ on the manifold
$M\times T^{n-k}$.
For $\gamma\in\Gamma$ let
$$
\gamma(u_1,\dots,u_n,y)=
(u_1\cos\pi(e_1,\gamma),\dots,u_n\cos\pi(e_n,\gamma),y+\gamma).
$$
Note that
$\cos\pi(e_j,\gamma)=\pm 1$.
We introduce a map
$$
 \varphi:M\times T^{n-k}\rightarrow \C^n,
$$
$$
 \varphi (u_1,\dots,u_n,y)=(u_1e^{\pi i(e_1,y)},\dots,u_ne^{\pi
 i(e_n,y)}).
$$
Let denote by $e$ a vector
$e_1+\dots+e_n$.
We assume $\C^n$ is endowed by the Euclidian metric and the
symplectic form
$$
\omega=dx_1\wedge dy_1+\dots+dx_n\wedge
dy_n.
$$
The following theorem is valid:

\begin{theorem}
Group $\Gamma$ acts on $M\times T^{n-k}$ freely.
Map $\varphi$ gives an immersion
$$
\psi_1:M_1=M\times T^{n-k}/\Gamma\rightarrow\C^n.
$$

Immersion $\psi_1$ is Lagrangian $H$-minimal.
If $e=0$, then $\psi_1$ is a Lagrangian minimal immersion.
\end{theorem}

Below we will show that $\psi_1$ is a bijective map
between everywhere open dense subsets in $M_1$ and
$\psi_1(M_1)$ which are determined by inequalities
$u_j\ne 0, j=1,\dots,n$.  Thus, selfintersections of $\psi_1(M_1)$
can appear only at points for which $u_j=0$ for some $j$.

For $M_1$ to be closed it is sufficient, for example, to require
that one of equations (2) gives an ellipsoid. The
condition $e=0$ does not allow $M_1$ to be closed.

The notion of $H$-minimality was introduced by Oh [1].
He also proved that the torus
$$
S^1(r_1)\times\dots\times S^1(r_n)\subset\C^n,
$$
where $S^1(r_j)$
is a circle of radius $r_j$ is an $H$-minimal Lagrangian
submanifold in $\C^n$. Other examples of $H$-minimal
Lagrangian tori in $\C^2$ have been construc\-ted by Castro and
Urbano [2] and Helein and Romon [3]. In [4] Helein and Romon
constructed an immersed (with selfintersection)
$H$-minimal Lagran\-gian Klein's bottle in $\C^2$. As was shown by
Nemirovski [5], there exist no Lagrangian Klein's
bottles embedded in $\C^2$.

Let's note that in $\C^n$ there exist no closed minimal
submanifolds.

We also note that Oh's example follows from the fact that the
circle $S^1(r)$ in $\C$ is an $H$-minimal Lagrangian submanifold,
and from the following lemma.

\begin{lemma}
Let
$P=P_1\times P_2$
and
$\omega=\pi_1^*\omega_1+\pi_2^*\omega_2$,
where
$(P_i,\omega_i)$
is K\"ahler manifold and
$\pi_i:P\rightarrow P_i$
is a proejection, $i=1,2$.
Let
$L_i\subset P_i$
be a $H$-minimal Lagrangian submanifold, $i=1,2$.
Then the submanifold
$L=L_1\times L_2 \subset P$
is $H$-minimal Lagrangian.
\end{lemma}

Lemma 1 is also valid for immersions.

Let's consider the case of $M$ being a cone, i.e. $d_j=0$ in
equation (2). Then to this cone there corresponds an
immersion $\psi_2$ into the complex projective space
$\C P^{n-1}$ of a $(n-1)$-dimensional manifold $M_2$ which
results from factorization of manifold
$(M-\{0\})\times T^{n-k}/\Gamma$ by action of
$\R^*$
$$ \psi_2(u_1,\dots,u_n,y)=(u_1e^{\pi
i(e_1,y)}:\dots:u_ne^{\pi i(e_n,y)}).
$$
We suppose that $\C P^{n-1}$ is endowed by the Fubini -- Studi
metric and as a symplec\-tic form we consider the K\"ahler form
$\Omega$ of this metric.

\begin{theorem}
Map $\psi_2$
is an $H$-minimal Lagrangian immersion.

If $e=0$, then $\psi_2$ is a minimal Lagrangian
immersion.
\end{theorem}

It is now evident how to construct the closed manifold $M_2$.
For this it is sufficient, for example, that in any of
the equations (2) all but one coefficient be positive.

Helein and Romon established a correspondence
between $H$-minimal Lagrangian cones in $\C^3$ and
$H$-minimal Lagrangian surfaces in $\C P^2$, but did not
propose explicit examples [6]. Castro and Urbano in [7]
constructed examples of Lagrangian minimal tori in $\C P^2$.

From Lemma 1 follows, that, using examples of $H$-minimal
Lagrangian immersions in $\C^n$ and $\C P^n$, one can construct
$H$-minimal Lagrangian immer\-sions in
$\C^n\times\C P^{n_1}\times\dots\times \C P^{n_k}.$

Lemma 1 is proven in paragraph 2, Theorems 1 and 2 are
proven in paragraphs 3 and 4. We will give some
examples for Theorems 1 and 2 at the end of paragraphs 3 and 4,
respectively.

The author thank Ya. V. Bazaikin and I. A.
Taimanov for valuable conversations.

\section{Proof of Lemma 1}

Let $H$ denote
a vector field of mean curvature along $L$.
Oh proved the following assertion:

{\sl A Lagrangian submanifold
$L\subset P$
is $H$-minimal
if and only if
$$
\delta\omega_H=0,
$$
where
$\delta$
is the Hodge dual operator to $d$.}

Actually, by formula of volume variation along
field $W$ we have
$$
\frac{d}{d t}{\rm vol}(L_t)=\int_L<H,W>d{\rm vol}=
\int_L<\omega_H, \omega_W>d{\rm vol},
$$
where $L_t$ is the deformation of $L$ along field $W$. For every
function $f$ on $L$ we can select a hamiltonian field $W$
such that $\omega_W=df$, consequently, this assertion follows
from the arbitrariness of $f$ (for details see [1]).

{\sc Proof of Lemma 1.} $L$ is Lagrangian because
$$
\omega(V,V')=\omega(\pi_1^*V_1+\pi_2^*V_2,\pi_1^*V'_1+\pi_2^*V'_2)=
$$
$$
\pi_1^*\omega_1(V_1,V'_1)+ \pi_2^*\omega_2(V_2,V'_2)=0,
$$
where
$V$
and
$V'$
are vectors tangent to
$L$,
$V_i$
and
$V'_i$
are uniquely definеd vectors to
$L_i$ such that $V$ and $V'$
decompose into the sum of their lifting on
$ L,\ V=\pi_1^*V_1+\pi_2^*V_2,
V'=\pi_1^*V'_1+\pi_2^*V'_2. $

We show that
$\delta
\pi_1^*\alpha=\pi_1^*\delta\alpha,$
where
$\delta=-*d*,\ *$
is Hodge's star,
$\alpha$ is a one-form on $P_1$.
Let
$(x)=(x_1,\dots,x_n)$ be coordinates in the neighbourhood of
$p_1\in P_1$, and let
$(y)=(y_1,\dots,y_k)$ be coordinates in the neighbourhood of
$p_2\in P_2.$
Because
$\delta$ and $\pi^*_1$ are linear it is sufficient to
consider the case  of
$\alpha=f(x)dx_1$.
Operator $*$ on $P_1$ has the form
$$
*\frac{dx_1}{g_1(x)}=\frac{1}{g_2(x)}dx_2\wedge\dots \wedge
dx_n,\ *\frac{1}{g_1(x)g_2(x)}dx_1\wedge\dots \wedge dx_n=1,
$$
where $g_1(x)$ and $g_2(x)$ are the norms of the forms
$dx_1$ and $dx_2\wedge\dots\wedge dx_n.$
Analogically, star $*$ on $P$
has the form
$$
*\frac{dx_1}{g_1(x)}=\frac{1}{g_2(x)g(y)}dx_2\wedge\dots \wedge
dx_n\wedge dy_1\wedge\dots\wedge dy_k,
$$
$$
 *\frac{1}{g_1(x)g_2(x)g(y)}dx_1\wedge\dots \wedge dx_n \wedge
dy_1\wedge\dots\wedge dy_k=1,
$$
where $g(y)$ is the norm of form
$dy_1\wedge\dots\wedge dy_k$.
Then
$$
\pi_1^*\delta\alpha=-\pi_1^**d*\alpha =
-\pi_1^**d\left(\frac{f(x)g_1(x)}{g_2(x)}dx_2\wedge\dots \wedge
dx_n\right)
$$
$$
=-\pi_1^**\left(\partial_{x_1}\left(\frac{f(x)g_1(x)}{g_2(x)}\right)dx_1\wedge\dots
\wedge dx_n\right)
=-\partial_{x_1}\left(\frac{f(x)g_1(x)}{g_2(x)}\right)g_1(x)g_2(x).
$$
On the other hand
$$
\delta\pi_1^*\alpha=-*d\left(\frac{f(x)g_1(x)}{g_2(x)g(y)}dx_2\wedge\dots
\wedge dx_n\wedge dy_1\wedge\dots\wedge dy_k\right)
$$
$$
=-*\left(\partial_{x_1}\left(\frac{f(x)g_1(x)}{g_2(x)g(y)}\right)dx_1\wedge\dots
\wedge dx_n\wedge dy_1\wedge\dots\wedge dy_k \right)
$$
$$
=-\partial_{x_1}\left(\frac{f(x)g_1(x)}{g_2(x)}\right)g_1(x)g_2(x).
$$
Analogically one can show that
$\delta
\pi_2^*\alpha=\pi_2^*\delta\alpha.
$
We note that by the definition of mean curvature it is not
difficult to obtain equality
$H=\pi_1^*H_1+\pi_2^*H_2$, where
$H_i$ is a vector field of mean curvature along $L_i$.
Then $H$-minimality of $L$
follows from equalities
$$
\delta\omega_H=\delta\omega_{\pi_1^*H_1+\pi_2^*H_2}=
\delta\pi_1^*{\omega_1}_{H_1}+\delta\pi_2^*{\omega_2}_{H_2}
$$
$$
=\pi_1^*\delta{\omega_1}_{H_1}+\pi_2^*\delta{\omega_2}_{H_2}=0.
$$
Thus, Lemma 1 is proven.

\section{Immersions in $\C^n$ (Theorem 1)}

There is a criterion of $H$-minimality in terms of
the Lagrangian angle --- a function defined on $L$:

{\sl  An immersion $\psi$ of manifold $L$ is
$H$-minimal if and only if the Lagrangian angle  is a
harmonic function on $L$.}

We recall the Wolfson construction [8] of building a
Lagrangian angle. We choose an orientation on $L$.
Let $K$ be a canonical linear bundle over $P$ with the
connection $\nabla$, induced by the Levi--Chivita connection.
Sections of $K$ are holomorphic $n$-forms on $P$.
We suppose that there exists some section $\sigma$ of
bundle $\psi^*K$ parallel on $L$ in induced connection, i.e.
$\psi^*\nabla \sigma=0$.  Let the norm of
$\sigma$ in each point be equal to 1.
In each point $x\in L$ we choose an orthonormal
tangent frame $\xi$ concordant with the orientation of $L$.
Then we obtain some function $\beta(x)$ on $L$,
defined by equality
$$
\sigma(\xi)=e^{i\beta(x)}.
$$
Function $\beta$ is called Lagrangian angle on $L$.
The value of $\beta$ is independent of the choice of $\xi$.
For another choice of $\sigma$ function $\beta$ can change
only by a constant, consequently, the one-form $d\beta$
is correctly defined on $L$ (for details see [8]). In the
general case function $\beta(x)$ is multivalued. By moving on a
cycle $\beta(x)$ can change its value by
$2k, k\in\Z$.
Wolfson [8] showed that $d\beta=\omega_H$.
Consequently, immersion $\psi$ is $H$-minimal if and only if the
Lagrangian angle is a harmonic function on $L$, i.e.
$\Delta\beta=\delta d\beta=0.$

With the help of the Lagrangian angle  one can calculate the
vector  of mean curvature of $\psi(L)$. From the
lagrangianality of $\psi$ and equalities
$$
(J_PH,V)_P=\omega(H,V)=d\beta(\psi^*(V))=(\psi_*({\rm
grad}\beta),V)_P,
$$
where $(.,.)_P$ is a Riemannian metric on
$P$, $J_P$ is a complex structure on $P$, $V$ is a tangent
vector to $\psi(L)$, we obtain
$$
H=-J_P\psi_*({\rm grad} \beta) \eqno{(3)}
$$
(in points of selfintersection of $\psi(L)$ there arise some
vectors $H$).

For Lagrangian immersions in $\C^n$, the Lagrangian angle can be
calculated by means of form
$$
dz_1\wedge\dots \wedge dz_n,
$$
since it is parallel on
$\C^n$, and, consequently,
$$
 \psi^*\nabla(\psi^*(dz_1\wedge\dots\wedge
dz_n))=0.
$$

The following lemma is necessary.

 Let there be given a block
diagonal metric of the form
$$
ds^2=\sum_{i,j=1}^kg_{ij}(x)dx_idx_j+
\sum_{i,j=1}^{n-k}\widetilde{g}_{ij}(x)dy_idy_j, \eqno{(4)}
$$
where $(x)=(x_1,\dots,x_k)$ and $(y) = (y_1,\dots,y_{n-k})$.

\begin{lemma}
If function $\alpha$ is linear on $y$:
$$
\alpha=\alpha_1y_1+\dots+\alpha_{n-k}y_{n-k}, \ \ \ \alpha_j\in\R,
$$
then it is harmonic in metric $ds^2$.
\end{lemma}

{\sc Proof.} We denote
${\rm det}g_{ij},{\rm det}\widetilde{g}_{ij}$ and the
components of the matrices $(g_{ij})^{-1}$ and
$(\widetilde{g}_{ij})^{-1}$ by  $g,\widetilde{g},g^{ij}$ and
$\widetilde{g}^{ij}$, respectively.
Then
$$
\Delta\alpha=\frac{1}{\sqrt{g\widetilde{g}}}
\sum_k\partial_{x_k}(\sum_ig^{ik}\sqrt{g\widetilde{g}}
\partial_{x_i}\alpha)+
\frac{1}{\sqrt{g\widetilde{g}}}
\sum_k\partial_{y_k}(\sum_i\widetilde{g}^{ik}\sqrt{g\widetilde{g}}
\partial_{y_i}\alpha)=0.
$$
Lemma 2 is proven.

By
$$
<\xi,\eta>=\sum_{i=1}^{n+1}\xi_i\overline{\eta}_i=(\xi,\eta)-
i\omega(\xi,\eta)
$$
we denote  the Hermitian product of vectors $\xi$ and $\eta$ in
$\C^n.$

{\sc Proof of Theorem 1.} Let
$\widetilde{F}=(f_1(x),\dots,f_n(x))$
be a map on $M,$ where $x\in U^k\subset\R^k.$
We consider a map on
$\psi_1(M_1)$
$$
 F=(f_1(x)e^{\pi
 i(e_1,y)},\dots,f_n(x)e^{\pi i(e_n,y)}):  U^k\times
 V^{n-k}\rightarrow \psi_1(M_1)\subset\C^n,
$$
where
$
  y\in V^{n-k}\subset\R^{n-k}.
$
 From equality (2) we obtain
$$
 \langle\partial_{y_j}F,\partial_{x_s}F\rangle=
 \pi ie_{1j}f_1\partial_{x_s}f_1+\dots+\pi ie_{nj}f_n\partial_{x_s}f_n=0,
 \eqno{(5)}
$$
consequently, the induced metric on
$M_1$
in local coordinates has the form (4).  Further, from (5),
$
\langle\partial_{x_s}F,\partial_{x_j}F\rangle\in\R
$
and
$
\langle\partial_{y_s}F,\partial_{y_j}F\rangle\in\R
$
we obtain equalities
$$
 \omega(\partial_{y_j}F,\partial_{x_s}F)=
 \omega(\partial_{y_j}F,\partial_{y_s}F)=
 \omega(\partial_{x_s}F,\partial_{x_j}F)=0,
$$
which means that the embedding $\psi_1$ is Lagrangian.

If manifold $M$ is smooth, then the vectors
$\partial x_{j}F(x)$
are linearly indepen\-dent, because
$$ <\partial x_jF(x),\partial x_kF(x)>=(\partial
x_j\widetilde{F}(x),\partial x_k\widetilde{F}(x)).
$$
Normal vectors to the hypersurfaces given by equations (2) have
the form
$$
 n_j=(e_{1j}u_1,\dots,e_{nj}u_n).
$$
If point
$(u_1,\dots,u_n)\in M$ is smooth, then $n_j$ are linearly
independent, conse\-quently, the matrix consisting of
components of vectors $n_j$ has rank  $n-k$.
By multiplying the $k$-th column by
$\pi ie^{\pi i(e_k,y)}$
its rank (over $\R$) does not change. Therefore the vectors
$$
 \partial_{y_j}F=(\pi ie_{1j}f_1e^{\pi i(e_1,y)},\dots,
 \pi ie_{nj}f_ne^{\pi i(e_n,y)})
$$
are linearly independent. From (5) follows that map $F$ is smooth
of maximal rank.

From (5) also follows that in tangent space of each point
one can choose an orthonormal basis of form
$$
 p_j=(p_{j1}(x)e^{\pi i(e_1,y)},\dots,p_{jn}(x)e^{\pi
 i(e_n,y)}), j=1,\dots,k,
$$
$$
 q_s=(\pi iq_{s1}(x)e^{\pi
 i(e_1,y)},\dots,\pi iq_{sn}(x)e^{\pi i(e_n,y)}), s=1,\dots,n-k.
$$
On $M_1$ one can choose an orientation such that the
Lagrangian angle has the form
$$
 \beta=-i\log (dz_1\wedge\dots
 \wedge dz_n (p,q))=-i\log i^{n-k}e^{i(e_1+\dots+e_n,y)},
$$
where
$(p,q)$ is a matrix (up to interchanging two lines)
consisting of compo\-nents of vectors $p_j$ and $q_s$. From Lemma 2
follows $\Delta\beta=0$, which proves the $H$-minimality of
$\psi_1$.  If $e=0$, then $\beta=const$ and from formula (3)
follows the minimality of $\psi_1$.

By definition of $\Gamma$ follows that if
$\gamma\in\Gamma, \gamma\ne 0$, then there is a vector
$e_j$ such that $(\gamma,e_j)$ is an odd number.
Indeed, if $(\gamma,e_s)\in 2\Z$ for all $s$,
then by definition $\frac{\gamma}{2}\in\Lambda^*$,
therefore $\gamma\in 2\Lambda^*$.  Consequently,
$e^{\pi i(e_j,y+\gamma)}=-e^{\pi i(e_j,y)}$, and thus the
action of group $\Gamma$ on $M\times T^{n-k}$ is free.

We show that selfintersections of
$\psi_1(M_1)$
can arise only in points for which
$u_j=0$ for some $j$.
Suppose that
$$ \varphi(u_1,\dots,u_n,y)=\varphi(u'_1,\dots,u'_n,y') $$
and
$u_j\ne 0$.
Then
$$ |u_1|=|u'_1|,\dots,|u_n|=|u'_n|, $$ $$
 e^{\pi i(e_1,y-y')}=\pm 1, \dots,
e^{\pi i(e_n,y-y')}=\pm 1,
$$
i.e.
$$
(e_1,y-y')\in\Z,\dots,(e_n,y-y')\in\Z,
$$
consequently, $\gamma=y-y'\in\Gamma$,
$$
u_1=u'_1\cos\pi(e_1,\gamma),\dots,u_n=u'_n\cos\pi(e_n,\gamma),
$$
$$
\gamma(u'_1,\dots,u'_n,y')=(u_1,\dots,u_n,y).
$$
Thus $\psi_1$ is an immersion of manifold $M_1$ in
$\C^n$. Theorem 1 is proven.

Example 1. Let $M$ be an ellipse
$$
u_1^2+2u_2^2=1.
$$
In this case
$$
\Lambda=\Z,\ \Lambda^*=\Z,\ \Gamma=\Z_2.
$$
The non-zero element $\gamma$ from $\Gamma$
acts on $M\times S^1$ in the following way:
$$
 \gamma(u_1,u_2,q)=(-u_1,u_2,-q).
$$
Up to a homeomorphism we can suppose that
$M$ is a unit circle in $\C^1$, then $\gamma$ acts on
torus $T=S^1\times S^1$ in the following way
$$
\gamma(q_1,q_2)=(\bar{q}_1,-q_2).
$$
We consider a universal covering of torus
$T$ by plane $\R^2$
$$
(x,y)\rightarrow (e^{2\pi ix},e^{2\pi iy}).
$$
An action $\Z_2$
on $T$ induces an action $Z_2$ on $\R^2$,
point $(x,y)$ moves by action of the non-zero element in
$(-x,y+\frac{1}{2})$. Thus
manifold $T/\Z_2$ is homeomorphic to
$\R^2/\{\Z\oplus\Z\oplus\Z_2\}$.
The rectangle with the vertices
$A(0,0),B(0,\frac{1}{2}),C(1,0)$, $D(1,1)$ serves as a fundamental
domain of action of group $\Z\oplus\Z\oplus\Z_2$ on $\R^2$.
Thereby vector $AB$ must be identified with
$CD$, and vector $AC$ --- with $BD$, but with a change of
orientation. Consequently,
$M\times S^1/\Gamma$ is Klein's bottle ${\cal K}$.  By mapping
$\varphi$ circles $\{ (0,\frac{1}{\sqrt{2}},y)\}$
and $\{ (0,-\frac{1}{\sqrt{2}},y+\frac{1}{2})\}$ of ${\cal K}$
are indentical.

Example 2. We consider a more general situation. Let
$M$ be a $(n-1)$-dimensional ellipsoid
$$
m_1u_1^2+\dots+m_nu_n^2=1,
$$
where $m_j$ are natural numbers.
The non-zero element $\gamma\in\Gamma=\Z_2$ acts on
$M\times S^1$ in the following way:
$$
\gamma(u,q)=(\tau(u),-q),
$$
where $\tau:M\rightarrow M$ is some involution.
We cut $S^1$ in two halves. Up to a homeomorphism
we can suppose that one of these parts is the interval $[0,1]$.
Then $M\times S^1/\Gamma$ is obtained from a cylinder
$M\times [0,1]$ by identification of points on
the boundary, namely, the point of the form
$(u,0)$ must be identified with
$(\tau (u),1)$.  If $\tau$ preserves the orientation
of $M$, then $M\times S^1/\Gamma$ is diffeomorphic to
$S^{n-1}\times S^1$, and if it does not, then $M\times S^1/\Gamma$ is the
generalized Klein's bottle ${\cal K}^n$.  For example,
in the case when $M$ is a sphere, the involution $\tau$
identifies point $u$ with $-u$. If the sphere is
evendimensional, then $\tau$ does not preserve its orientation,
but if the sphere is odddimensional, then $\tau$
does preserve its orientation. Thereby map $\varphi$
induces the embedding of manifold $M\times S^1/\Gamma$.

Example 3.  Let $M$ be an intersection of sphere
$$
 u_1^2+\dots+u_n^2=1
$$
and cone
$$
 u_1^2+\dots+u_{n-1}^2=u_n^2,
$$
$M$ is a unification of two non-intersecting spheres
 ($u_n>0$ and $u_n<0$).
In this case
$$
 \Lambda=e_1\Z\oplus e_3\Z,\ e_1=(1,1),\
 e_3=(1,-1),
$$
$$ \Lambda^*=\gamma_1\Z\oplus \gamma_2\Z,
  \gamma_1=\left(\frac{1}{2},\frac{1}{2}\right),
  \gamma_2=\left(\frac{1}{2},-\frac{1}{2}\right), \ \Gamma=\Z_2\oplus\Z_2.
$$
The elements
 $\gamma_1$ and $\gamma_2$ act on
 $M\times   T^2$
 in the following way
 $$
  \gamma_1(u_1,\dots,u_n,y)=\left(u_1\cos\pi,\dots,u_{n-1}\cos\pi,u_n\cos
  0,y+\gamma_1\right)
 $$
 $$
  =\left(-u_1,\dots,-u_{n-1},u_n,y+\gamma_1\right),
 $$
 $$
 \gamma_2(u_1,\dots,u_n,y)
  =\left(u_1,\dots,u_{n-1},-u_n,y+\gamma_2\right).
 $$
The element $\gamma_2$ glues two non-intersecting components
$M\times T^2$ ($u_n>0$ and $u_n<0$). Up to a homeomorphism
one can suppose that the component of connectedness
$M\times T^2$  is $S^{n-2}\times S^1\times S^1,$
then an action $\Z_2$ on it looks the following way
$$
\gamma  (q_1,q_2,q_3)=(-q_1,-q_2,q_3).
$$
Consequently,
$M\times T^2/\Gamma$ is diffeomorphic to
${\cal K}^{n-1} \times S^1$ if $n$ is an even number
and $S^{n-2}\times S^1\times S^1$ if $n$ is an odd
number (see  example 2).  We prove that $\psi_1(M_1)$
has no selfintersections.  We suppose that
$$
\varphi(u_1,\dots,u_n,y)=\varphi(u'_1,\dots,u'_n,y').
$$
It is clear that
$u_n\ne 0$ and at least one of the components
$u_1,\dots,u_{n-1}$ is non-zero. Let $u_j\ne 0$.
Then
$$
e^{\pi i(e_j,y-y')}=\pm 1, e^{\pi i(e_n,y-y')}=\pm 1,
$$
i.e.
$$
 (e_j,y-y')\in\Z,\ (e_n,y-y')\in\Z,
$$
Since
$e_1=\dots=e_{n-1}$, then $\gamma=y-y'\in\Gamma$,
$$
u_1=u'_1\cos\pi(e_1,\gamma),\dots,u_n=u'_n\cos\pi(e_n,\gamma).
$$
Consequently,
$$
\gamma(u'_1,\dots,u'_n,y')=(u_1,\dots,u_n,y).
$$
Example 4. Let $M$ be an intersection of ellipsoid
$$
  u_1^2+2u_2^2+u_3^2=1
$$
and cone
$$
 u_1^2+2u_2^2=u_3^2.
$$
As in the previous example an action $\Gamma$ on
$M\times T^2$ is reduced to an action $\Z_2$  on torus
$S^1\times S^1\times S^1\subset \C^3$
$$
\gamma(q_1,q_2,q_3)=(\bar{q_1},-q_2,q_3).
$$
Consequently,
$M\times T^2/\Gamma$ is diffeomorphic to
${\cal K}\times S^1$.

\section{Immersions in $\C P^n$ (Theorem 2)}

{\sc Proof of Theorem 2.}
We consider the Hopf bundle
$$
h:S^{2n-1}\rightarrow\C P^{n-1}.
$$
We suppose that $S^{2n-1}$ is a unit sphere
in $\C^n$.  $h$ maps every circle $S^1\subset S^{2n-1}$
given by an intersection of $S^{2n-1}$ with the complex line
$l\subset\C^n,$ $0\in l$ into a point of the
projective space corresponding to $l$.

Let $\widetilde{L}$ be a simple connected
$(n-1)$-dimensional submanifold in $S^{2n-1}$
such that for any point $p\in\widetilde{L}$
the linear hull of radius-vector
$p$ and the tangent space to
$\widetilde{L}$ in point $p$ is a Lagrangian
$n$-dimensional space in $\C^n$. We shall suppose that
$\widetilde{L}$ is a sufficiently small neigbourhood of
point $p$.
By $L$ we denote $h(\widetilde{L})\subset \C P^{n-1}.$

\begin{lemma}
Submanifold  $L\subset \C P^{n-1}$
is Lagrangian. Map $h:\widetilde{L}\rightarrow L$
is an isometry.
\end{lemma}

{\sc Proof.} In every point
$p\in \widetilde{L}$ we choose a tangent basis
$\xi_1,\dots,\xi_{n-1}$ to $\widetilde{L}$.
Then by condition of Lemma 3 we have equalities
$$ \omega(\xi_k,p)=(\xi_k,ip)=(\xi_k,p)=0.  $$
Consequently, every vector $\xi_k$
is orthogonal to circle
$S^1_p\subset S^{2n-1},$ which is an intersection of
$S^{2n-1}$ and a complex line passing through
$p$. Thus the tangent space to $\widetilde{L}$ in point
$p$ is orthogonal to a fiber of the Hopf bundle
including $p$, i.e.
manifold $\widetilde{L}$ is a horizontal
manifold by Hopf mapping $h$.  Consequently, since
$h$ is a Riemannian submersion,
$h|_{\widetilde{L}}$ is an isometry.
Equality
$$
\Omega(h_*(\xi_k),h_*(\xi_j))=\omega(\xi_k,\xi_j)=0
$$
holds, which means that $L$ is a Lagrangian submanifold
in $\C P^{n-1}$.  Lemma 3 is proven.

We define a $(n-1)$-form $\sigma$ on $L$. For any tangent
vectors $\eta_1,\dots,\eta_{n-1}$ to $L$ in point
$h(p)$, where $p\in\widetilde{L}$ there exist
uniquely defined tangent vectors $\xi_1,\dots,\xi_{n-1}$ to
$\widetilde{L}$ in point $p$ such that
$h_*(\xi_i)=\eta_i$.
We suppose
$$
\sigma(\eta_1,\dots,\eta_{n-1})=
\sigma(h_*(\xi_1),\dots,h_*(\xi_{n-1}))= dz_1\wedge\dots\wedge
dz_n(\xi_1,\dots,\xi_{n-1},p).
$$
We show that with the help of
$\sigma$ one can calculate the Lagrangian angle
$\beta_1$ on $L$.  For this it is necessary
to show that $\nabla |_L\sigma=0$, where $\nabla$
is a connection in a canonical bundle over $\C P^{n-1}$,
concordant with the Levi--Chivita connection.
Analogically one can define
$\sigma$ on the tangent vectors
$\eta_1,\dots,\eta_{n-1}$ in points of $L$ to $\C P^{n-1}$.
For this it is necessary to take their horizontal lifting
$\xi_1,\dots,\xi_{n-1}$ on $\widetilde{L}$.

\begin{lemma}
Form $\sigma$ is parallel on $L$
for connection $\nabla|_L$ .
\end{lemma}

{\sc Proof.} In order to prove parallelism of $\sigma$
it is sufficient to consider an arbitrary smooth path $s(t)$ in
$L$, to consider an arbitrary parallel tangent frame
$\eta(t)=(\eta_1(t),\dots,\eta_{n-1}(t))$ to $\C P^{n-1}$ along
$s$ and to show that $\sigma(\eta(t))$ is not dependent on
parameter $t$. Let
$\widetilde{s}\subset\widetilde{L}$ be a horizontal lifting of
$s$ and $(\xi_1(t),\dots,\xi_{n-1}(t))$ be a
horizontal lifting of frame $\eta(t)$ by $h$.
For a Levi--Chivita connection $D$ on $S^{2n-1}$ and
vector fields $X$ and $Y$ we have
$$
D_XY={\cal H}D_XY+{\cal V}D_XY,
$$
where ${\cal H}D_XY$ and ${\cal V}D_XY$
are horizontal and vertical components, respective\-ly.
Moreover, $h_*({\cal H}D_XY)=\widetilde D_{h_*(X)}h_*(Y)$, where
$\widetilde D$ is a Levi--Chivita connec\-tion on $\C P^{n-1}$.
Consequently, since field $\eta_i(t)$ is parallel along $s$ and
$\xi_i(t)$ is its horizontal lifting, field
$\frac{d\xi_i}{dt}$ has the form $f_i(t)\widetilde{s}(t)$, where
$f_i(t)$ is some complex-valued function, since the projection
of vector $\frac{d\xi_i}{dt}$ on $S^{2n-1}$ must not
have a horizontal component.  The following equality holds:
$$
\sigma(\eta(t))={\rm det} \left( \begin{array}{ccc} \xi_1^1(t) &
\dots & \xi_1^n(t)\\ \dots & \dots & \dots \\ \xi_{n-1}^1(t) &
\dots & \xi_{n-1}^n(t)\\ s^1(t) & \dots & s^n(t)
\end{array}\right),
$$
where
$(\xi_i^1,\dots,\xi_i^n)$ are complex coordinates of
vectors
$ \xi_i,\  (s^1(t),\dots,s^n(t)) $ are complex coordinates of
$s$.  Then
$$
 \frac{d}{d t}\sigma(\eta(t))={\rm det}
 \left( \begin{array}{ccc}
 \frac{d}{d t}\xi_1^1(t) & \dots & \frac{d}{d t}\xi_1^n(t)\\
\dots & \dots & \dots \\ \xi_{n-1}^1(t) & \dots &
 \xi_{n-1}^n(t)\\ s^1(t) & \dots & s^n(t)
 \end{array}\right)+\dots+
 $$
 $$ {\rm det} \left(
 \begin{array}{ccc}
 \xi_1^1(t) & \dots & \xi_1^n(t)\\
 \dots & \dots & \dots \\
 \frac{d}{d t}\xi_{n-1}^1(t) & \dots & \frac{d}{d
 t}\xi_{n-1}^n(t)\\
 s^1(t) & \dots & s^n(t)
 \end{array}\right)+
{\rm det}
 \left(
 \begin{array}{ccc}
 \xi_1^1(t) & \dots & \xi_1^n(t)\\
 \dots & \dots & \dots \\
 \xi_{n-1}^1(t) & \dots & \xi_{n-1}^n(t)\\
 \frac{d}{d t}s^1(t) & \dots & \frac{d}{d t}s^n(t)\\
 \end{array}\right)=0.
$$
The first $n-1$ determinants are equal to 0 by virtue of
linear dependency of vectors
$\frac{d}{d t}\xi_i$ and $s$.
Vectors $\eta_1,\dots,\eta_{n-1},\frac{ds}{dt}$ are
(complex) linearly dependent, consequently, their
horizontal liftings
$\xi_1,\dots,\xi_{n-1},\frac{d
\widetilde{s}}{d t}$ are also linearly depen\-dent,
therefore the last determinant is equal to zero, as well.
Lemma 4 is proven.

Let $p$ be some intersection point of cone
$\psi_1(M_1)\subset\C^n$ and  sphere $S^{2n-1}$. By $m$
we denote one of the points in the inverse image
$\psi_1^{-1}(p)$ (if $p$ is a selfintersection point of
$\psi_1(M_1)$, then there are several).
By $V$ we denote some small neighbourhood of $m$ and by
$\widetilde{L}$ the intersection $\psi_1(V)\cap S^{2n-1}$.
By Lemma 3 submanifold $h(\widetilde{L})\subset\C P^{n-1}$
is Lagrangian.
Consequently,  immersion $\psi_2$ is also Lagrangian,
because $\psi_2(M_2)$ can be covered by  neighbourhoods of type
$h(\widetilde{L})$. By construction
of form $\sigma$ we have
$$
 \beta(p)=\beta_1(h(p))
$$
(for simplicity we will identify $V$ and $\psi_1(V)$ and also
small neighbourhoods in $M_2$ with their images by  $\psi_2$).

By choosing map $F$ in the proof of Theorem 1 one can suppose
that one of the coordinates, say, $x_1$, is a coordinate along
straight lines that generate a cone, and
$x_2,\dots,x_k$ are coordinates on $\psi_1^{-1}(\widetilde{L})$.
Then the metric on $\widetilde{L}$,
and, as follows from Lemma 3, on $L$, in the neighbourhood of point
$h(p)$ has the form (4), where $x=(x_2,\dots,x_k)$
and the summation on coordinate $x_1$ are skipped.
By Lemma 2 function $\beta_1$ is harmonic on $M_2$, which
proves the $H$-minimality of $\psi_2$. From formula (3)
follows the remaining part of the theorem. Theorem 2 is proven.

Example 5. Let $k=1, d_j=0$ in (2). Then $\psi_2(M_2)$
is a $(n-1)$-dimensional Clifford torus in
$\C P^{n-1}$
$$ \psi_2(M_2)=\{(u_1e^{\pi i(e_1,y)}:\dots:u_ne^{\pi
i(e_n,y)})\}\subset {\C P}^{n-1}.
$$

Example 6. Let $M$ be a cone
$$
u_1^2+2u_2^2=3u_3^2.
$$
Then
$$
\Lambda=\Z,\ \Lambda^*=\Z,\ \Gamma=\Z_2.
$$
The non-zero element from $\Gamma$ maps point
$(u_1,u_2,u_3,q)\in M\times S^1$ into point
$(-u_1,u_2,-u_3,-q)$.  Consequently, $\psi_2$
is a minimal immersion of Klein's bottle
(see example 1).

Example 7. Let $M$ be given by equation
$$
m_1u_1^2+\dots+m_{n-1}u_{n-1}^2=m_nu_n^2,
$$
where $m_j$ are natural numbers. In this case the topological
type of $M_2$ will be the same as in example 2. If
$m_1+\dots+m_{n-1}=m_n$, then immersion $\psi_2$
is minimal. If $m_1=\dots = m_{n-1},$ then  immersion
$\psi_2$ is an embedding.

It is not difficult to construct examples in ${\C P}^n$
analogically to those in 3 and 4.

\vskip7mm

{\bf BIBLIOGRAPHY }

\vskip3mm

[1] Oh, Y. Volume minimization of Lagrangian submanifolds under
Hamil\-to\-nian deformations // Math. Z. 1993. V. 212. P. 175-192.

[2] Castro, I., and Urbano, F. Examples of unstable Hamiltonian-minimal
Lagrangian tori in $\C^2$ // Compositio Math. 1998. V. 111.
P. 1-14.

[3]  Helein, F., and Romon, P. Hamiltonian stationary Lagrangian
surfaces in $\C^2$ //
Comm. Anal. Geom. 2002. V. 10. P. 79--126.

[4] Helein, F., and Romon, P. Weierstrass representation of
Lagrangian surfaces in four-dimensional space using spinors and
quaternions // Commen\-tari Mathematici Helvetici. 2000. V. 75. P.
680-688.

[5] Nemirovski, S.
 Lefschetz pencils, Morse functions, and Lagrangian embeddings of the Klein bottle
 //  Izvestia Math. 2002. V. 66. N. 1. P. 151-164.

[6] Helein, F., and Romon, P. Hamiltonian stationary Lagrangian
surfaces in Hermitian symmetric spaces // In: Differential
geometry and Integrable Systems. Eds. M. Guest, R. Miyaoka, and
Y. Ohnita. Contemporary Mathe\-matics. V. 308. Amer. Math. Soc.,
Providence, 2002. P. 161-178.

[7] Castro, I., and Urbano, F. New examples of minimal Lagrangian
tori in the complex projective plane // Manuscripta Math. 1994.
V. 85. P. 265--281.

[8] Wolfson, J. Minimal Lagrangian diffeomorphisms and the
Monge-Ampere equation // J. Differential Geometry. 1997. V. 46.
P. 335-373.

\end{document}